\documentclass[2p]{article}
\usepackage{amssymb,amsmath,amsthm}
\usepackage{indentfirst}
\usepackage{exscale}
\usepackage{relsize}
\usepackage[numbers,sort&compress]{natbib}

\usepackage{geometry}
\usepackage{color}
\usepackage[colorlinks, bookmarksopen, bookmarksnumbered, citecolor=black, citecolor=blue ]{hyperref}
\newcommand{\R}{\mathbb{R}}
\theoremstyle{definition}

\allowdisplaybreaks
\theoremstyle{remark}

\numberwithin{equation}{section}
 \UseRawInputEncoding

\begin{document}
\title{\Large\bf{ Infinitely many solutions for an instantaneous and non-instantaneous fourth-order differential system with local assumptions}}
\date{}
 \author {Lijuan Kang$^{1}$,  \ Xingyong Zhang$^{1,2}$\footnote{Corresponding author, E-mail address: zhangxingyong1@163.com}, \ Cuiling Liu$^{1}$\\
{\footnotesize $^1$Faculty of Science, Kunming University of Science and Technology,}\\
 {\footnotesize Kunming, Yunnan, 650500, P.R. China.}\\
{\footnotesize $^{2}$Research Center for Mathematics and Interdisciplinary Sciences, Kunming University of Science and Technology,}\\
 {\footnotesize Kunming, Yunnan, 650500, P.R. China.}\\}
 \date{}
 \maketitle

 \begin{center}
 \begin{minipage}{15cm}
 \par
 \small  {\bf Abstract:}  We investigate a class of fourth-order differential systems with instantaneous and non-instantaneous impulses. Our technical approach is mainly based on a variant of Clark's theorem without the global assumptions. Under locally subquadratic growth conditions imposed on the nonlinear terms $f_i(t,u)$ and impulsive terms $I_i$, combined with perturbations governed by arbitrary continuous functions of small coefficient $\varepsilon$, we establish the existence of multiple small solutions. Specifically, the system exhibits infinitely many solutions in the case where $\varepsilon=0$.
 \end{minipage}
 \end{center}
  \allowdisplaybreaks
 \vskip2mm
 {\section{Introduction and main results}}
The study of impulse phenomena has evolved from instantaneous to non-instantaneous impulse models. In early research, Milman-Myshkis \cite{V.MilMan 1960} proposed the concept of instantaneous impulses, which describes abrupt changes (e.g., circuit tripping) occurring in negligible time. These models were widely applied in physics and engineering. However, subsequent studies revealed limitations in explaining natural phenomena like earthquakes or tsunamis, which persist over minutes. To address this, Hern$\acute{\mbox{a}}$ndez-O'Regan \cite{E. Hern¨¢ndez 2013} introduced a novel non-instantaneous impulse system, emphasizing prolonged durations and state-dependent dynamics of such impulses. This breakthrough spurred rapid advancements in the field over the past decade. Agarwal-O'Regan \cite{Agarwal R 2017} formally categorized impulse systems into instantaneous and non-instantaneous types, establishing a theoretical foundation for further research. More details on these two types are available in \cite{Y. Tian 2013, B. Ahmad 2009, Colao V 2016, Y.F. Wei 2019, Zhao 2020}. In the study of impulsive differential equations,  variational approach, fixed-point theory, topological degree theory, analytic semi-group and comparison method have been systematically developed as fundamental parts of an established methodological framework. We refer the readers to \cite{Heidarkhani S 2020, Bai L 2017, Heidarkhani S 2018, Colao V 2016, Pierri M 2013, Feckan M 2014, Qian D B 2005, Liu Z H 2013}.
\par
In recent years, the existence and multiplicity of solutions for impulsive differential equations incorporating both instantaneous and non-instantaneous impulses have attracted significant attention from the academic community. In \cite{Y. Tian 2019}, Tian-Zhang considered a class of differential equations with instantaneous and non-instantaneous impulses:
 \begin{equation}
 \label{eq0.1}
 \begin{cases}
 -u''(t)=f_i(t,u(t)),\; t\in(s_i,t_{i+1}],\; i=0,1,2,...,N,\\
 \Delta u'(t)=I_i(u(t_i)),\; i=1,2,...,N,\\
 u'(t)=u'(t_i^+),\; t\in(t_i,s_i],\; i=1,2,...,N,\\
 u'(s_i^+)=u'(s_i^-),\; i=1,2,...,N,\\
 u(0)=u(T)=0,
  \end{cases}
 \end{equation}
where $0=s_0<t_1<s_1<t_2<...<s_N<t_{N+1}=T$, $f_i\in C((s_i,t_{i+1}]\times \mathbb{R}, \mathbb{R})$, $I_i\in C(\mathbb{R},\mathbb{R})$, $\Delta u^{\prime}(t_i)=u^{\prime}(t_i^+)-u^{\prime}(t_i^-)$, the instantaneous impulses occur at the points $t_i$ and the non-instantaneous impulses continue on the intervals $(t_i,s_i]$.  They obtained that the problem (\ref{eq0.1}) admits at least one classical solution via variational method.
Based on \cite{Y. Tian 2019}, Zhou-Deng-Wang \cite{Zhou 2020} by using the critical point theory investigated the existence of solutions for fractional differential equations of $p$-Laplacian with instantaneous and non-instantaneous impulses. For further studies on the existence and multiplicity of solutions for impulsive differential equations, one can see references \cite{Li D P 2023, Zhang W 2020, Bonanno G 2014, Fan X L 2023, Qiao Y 2021}  and the references therein.

\par
Recently, in \cite{Xia 2023},  Xia-Zhang-Xie investigated the following fourth-order nonhomogeneous impulsive differential system with $\varepsilon=0$:
  \begin{equation}
 \label{eq1.1}
 \begin{cases}
 u^{(4)}(t)+2h(t)u^{\prime \prime \prime}(t)+\left(h^{2}(t)+h^{\prime}(t)+\beta e^{-\int_0^t h(\tau)d\tau}\right)u^{\prime \prime}(t)=f_i(t,u(t))+\varepsilon g_i(t,u(t)),\; t\in (s_i,t_{i+1}],\\
 i=0,1,2,...,N,\\
 u^{(4)}(t)+2h(t)u^{\prime \prime \prime}(t)+\left(h^{2}(t)+h^{\prime}(t)+\beta e^{-\int_0^t h(\tau)d\tau}\right)u^{\prime \prime}(t)=0,\; t\in (t_i,s_i],\; i=1,2,...,N,\\
 u^{\prime \prime}(t_i^+)=u^{\prime \prime}(t_i^-),\ s^{\prime \prime}(t_i^+)=s^{\prime \prime}(t_i^-),\ u^{\prime \prime \prime}(s_i^+)=u^{\prime \prime \prime}(s_i^-),\; i=1,2,...,N,\\
 \Delta u^{\prime \prime \prime}(t_i)=I_i(u(t_i)),\; i=1,2,...,N,\\
 u(0)=u^{\prime}(0)=u(T)=u^{\prime}(T)=0,
 \end{cases}
 \end{equation}
where $0=s_0<t_1<s_1<t_2<...<s_N<t_{N+1}=T$, $\Delta u^{\prime \prime \prime}(t_i)=u^{\prime \prime \prime}(t_i^+)-u^{\prime \prime \prime}(t_i^-)$, $u^{\prime \prime \prime}(t_i^\pm)=\lim_{t\to t_i^\pm}u^{\prime \prime \prime}(t)$ for all $i=1,2,...,N$, $h(t)\in C^1([0,T];\mathbb{R})$, $\beta\in\mathbb{R}$ is a constant, and $\varepsilon$ is a small parameter. In \cite{Xia 2023}, the authors studied the specific case where $\beta=0$ and $\varepsilon=0$. They assumed the nonlinear term $f_i(t,u)$ is odd and satisfies superquadratic growth conditions, while the impulse term $I_i(u)$ is odd and subquadratic. Using the Mountain Pass Theorem, they proved the existence of at least one non-trivial weak solution. Further, by applying the Symmetric Mountain Pass Theorem, they showed the system (\ref{eq1.1}) admits infinitely many solutions. However, these conclusions are strictly limited to the special case $\beta=\varepsilon=0$, leaving room for extensions to general cases where $\beta\neq 0$ or $\varepsilon\neq 0$.

\par
Motivated by the above fact, in this paper, we consider the multiplicity of solutions for system (\ref{eq1.1}) with the case that $\beta\not=0$ and $\varepsilon\neq 0$, and  assume that the nonlinear terms $f_i(t,u)$ and impulsive terms $I_i(u)$ are odd functions with  subquadratic growth conditions only near $u=0$ and $g_i(t,u)$ ia allowed to be any continuous function near $u=0$. Such problem has two difficulties: (1) the variational functional becomes undefined and lacks smoothness on the working place; (2) the perturbation term $\varepsilon g_i$ break the symmetry of the equation (\ref{eq1.1}). As a result, the classical Clark's Theorem (\cite{P.H. Rabinowitz}, Theorem 9.1) cannot be directly applied to equation (\ref{eq1.1}).
To overcome the first challenge, a modified Clark's Theorem from \cite[Theorem A. and Theorem 1.1]{Liu Z L 2015} can effectively solve this issue.
However, the presence of perturbation term $\varepsilon g_i$ renders the previous method ineffictive.
This challenging problem was addressed by Kajikiya in \cite{Kajikiya R 2013} using a new version of Clark's theorem and then together with cut-off technique.
Subsequently, the approach in \cite{Kajikiya R 2013} has been extended to various problems, including $p$-Kirchhoff type equations \cite{Liu J Y 2017}, semilinear Schr$\rm{\ddot{o}}$dinger equations \cite{Zhang Q Y 2024}, and quasilinear elliptic equations or systems \cite{Liu C L 2024}.
In this work, we will apply  Kajikiya's variant of Clark's theorem in \cite{Kajikiya R 2013}, which lacks local symmetry conditions,  to the impulsive differential system (\ref{eq1.1}) with the locally subquadratic growth for the nonlinear terms $f_i(t,u)$ and impulsive terms $I_i(u)$. This approach allows us to overcome the lack of symmetry induced by the perturbation term $\varepsilon g_i$ and establish the existence of multiple solutions.

\par
To clarify the main results, we need the following assumptions.
\begin{itemize}
\item[$(F_0)$]   $f_i, g_i\in C((s_i,t_{i+1}] \times [-\delta,\delta], \mathbb{R})$ for some $\delta>0$, $i=0,1,2,...,N$,  and $I_i\in C([-\delta,\delta],\mathbb{R})$ for all $i=1,2,...,N$;
\item[$(F_1)$] $f_i(t,u)=-f_i(t,-u), (t,u)\in (s_i,t_{i+1}]\times [-\delta,\delta]$, $i=0,1,2,...,N;$
 \item[$(F_2)$]  $F_i(t,u)\ge 0$, $uf_i(t,u)-2F_i(t,u)<0$ for all $0<\lvert u\rvert \leq\delta$ and $t\in (s_i,t_{i+1}], $
 where $F_i(t,u)=\int_0^u f_i(t,s)ds$, $i=0,1,...,N$;
\item[$(F_3)$]$\lim\limits_{|u| \to 0}\inf\limits_{t\in (s_i,t_{i+1}]}\dfrac {F_{i}(t,u)}{|u|^2}=+\infty$,  $i=0,1,2,...,N$;
\item[$(F_4)$]there exist constants $a_1>0$ and $1< \theta_1<2$  such that
$$\lvert f_i(t,u)\rvert\leq a_1\lvert u\rvert^{\theta_1-1},\ (t,u)\in (s_i,t_{i+1}]\times [-\delta,\delta],\ i=0,1,2,...,N;$$
\item[$(F_5)$]there exists a  constant $a_*$ with $0<a_*<\frac{a_1}{\theta_1}$ such that
$$F_i(t,u)\ge a_*|u|^{\theta_1}\;\text{for all }\;t\in(s_i,t_{i+1}]\;\text{and}\;u\in [-\delta,\delta],\ i=0,1,2,...,N;$$

\item[$(I_1)$]$I_i(u)=-I_i(-u), \mbox{for all } u\in[-\delta,\delta],\ i=1,2,...,N;$
\item[$(I_2)$]there exist constants  $a_2>0$ and $1<\theta_2<2$ such that
$$\lvert I_i(u)\rvert\leq a_2(\lvert u\rvert^{\theta_2-1}+1),\;\text{for\;all}\; u\in [-\delta,\delta],\ i=1,2,...,N;$$
\item[$(I_3)$] there holds\\
\hspace*{1.3cm} (i)\  $  I_i(u)\geq 0$, for all  $u\in [0, \delta]$, $i=1,2,...,N$,  \\
\hspace*{1.3cm} (ii)\
$ 2\sum\limits_{i=1}^N\int_0^{u}e^{H(t_i)}I_i(t){\rm d}t\geq  \sum\limits_{i=1}^Ne^{H(t_i)}I_i(u)u, \mbox{ for all } u\in [-\delta, \delta],\ i=1,2,...,N,$\\
where $H(t)=\int_0^t h(\tau)d\tau$;

\item[$(G)$]$g_i:(s_i,t_{i+1}] \times [-\delta,\delta]\to \mathbb{R}$ is a continuous function  and there exists a function $b\in C([0,\delta],\R^+)$ such that
    $$
    |g_i(t,u)|\le b(|u|),\ \mbox{ for all }t\in (s_i,t_{i+1}],\  i=0,1,2,...,N;
    $$

\item[$(H)$]the following inequality holds:
\begin{equation*}
%\label{eq100a}
\frac{2e^{H_0}}{T^2}>\beta,
\end{equation*}
where $H_0=\min_{t\in [0,T]}{H(t)}$.
\end{itemize}
 \par

\vskip2mm
\noindent
{\bf Theorem 1.1.} {\it Suppose that assumptions $(F_0)$-$(F_5)$, $(I_1)$-$(I_3)$, $(G)$ and $(H)$ hold. Then for any $k\in \mathbb{N}$, there exists a constant $\varepsilon(k)>0$ such that if $|\varepsilon|\leq \varepsilon(k)$, system (\ref{eq1.1}) possesses at least $k$ distinct solutions whose $L^{\infty}$-norms are less than $\frac{\delta}{2}$. }

\vskip2mm
The following corollary is an immediate consequence of Theorem 1.1.
\vskip2mm
\noindent
{\bf Corollary 1.1.} {\it Assume that $(F_0)$-$(F_5)$, $(I_1)$-$(I_3)$  and $(H)$ hold. Then system (\ref{eq1.1}) with $\varepsilon=0$ possesses infinitely many distinct solutions whose $L^{\infty}$-norms are less than $\frac{\delta}{2}$.}

\vskip2mm
\noindent
{\bf Remark 1.1.} Unlike \cite{Xia 2023}, in this paper we consider the case where  $\varepsilon$ is allowed to be unequal to $0$. Thus, when the function $g_i(t,u)$ violates the odd function property $(g_i(t,-u)\neq -g_i(t,u))$, this implies that the local symmetry of system (\ref{eq1.1}) is broken. Such a case is not considered in \cite{Xia 2023}, and our conditions are still different from those in \cite{Xia 2023} even if $\varepsilon= 0$, because all of our assumptions on the nonlinear terms $f_i(t,u)$ and the impulsive terms $I_i(t,u)$ are defined only near the origin with respect to $u$ and $f_i(t,u)$ admit subquadratic growth at origin, whereas Theorem 3.7 of \cite{Xia 2023} showed that the unperturbed problem $(\varepsilon=0)$ has infinitely many weak solutions when $f_i(t,u)$ and $I_i(t,u)$ are  defined  globally for all $u\in\mathbb{R}$ and $f_i(t,u)$ admit superquadratic growth at infinity. Crucially, unlike \cite{Xia 2023}, which imposes the restrictive condition $\beta=0$, our approach permits $\beta<\frac{2e^{H_0}}{T^2}$.

\vskip2mm
\par
This paper is organized as follows: in Section 2, we define the working space and provide the necessary lemmas for the subsequent proofs. In Section 3, we present the multiplicity of non-trivial week solutions for (\ref{eq1.1}) with $\beta<\frac{2e^{H_0}}{T^2}$ and complete the proof of Theorem 1.1. In section 4 we give an example that illustrate our results.

\section{Preliminaries}\label{section 2}
  \setcounter{equation}{0}
Consider the space $E=H_0^2(0,T)$, where
\begin{align*}
&E=\{u:(0,T)\to \mathbb{R} |\ u\;\text{and}\ u^{\prime} \; \text{are\;absolutely\;continuous\;on}\;(0,T),\ u^{\prime\prime}\in L^2(0,T),\notag\\
&u(0)=u(T)=u^{\prime}(0)=u^{\prime}(T)=0\}.
\end{align*}
$E$ is a Hilbert space with the inner product
$$
\langle u,v\rangle:=\int_0^{T}e^{H(t)}u^{\prime\prime}(t)v^{\prime\prime}(t){\rm d}t
$$
and the norm
$$
\|u\|_E:=\left(\int_0^Te^{H(t)}|u^{\prime\prime}(t)|^2{\rm d}t\right)^{1\over2},
$$
which is equivalent to the norm
$$
\|u\|:=\left(\int_0^T|u^{\prime\prime}(t)|^2{\rm d}t\right)^{1\over2}.
$$
Define
$$
C^1([0,T])=\{u:[0,T]\to \mathbb{R} |\ u\;\text{is\;continuously\;differentiable on }[0,T]\}
$$
with the norm
$$\|u\|_{C^1}=\max\left\{\|u\|_\infty, \|u^{\prime}\|_\infty\right\}.$$\\
For each $\varepsilon\in [0,1]$, we define the functional $J_\varepsilon(u):E\to \mathbb{R}$ as
\begin{eqnarray*}
J_\varepsilon(u)&=&\frac{1}{2}\|u\|_E^2-\frac{\beta}{2}\int_0^T|u^{\prime}(t)|^2{\rm d}t-\sum\limits_{i=0}^N\int_{s_i}^{t_{i+1}}e^{H(t)}F_i(t,u(t)){\rm d}t\notag\\
&& -\sum\limits_{i=1}^N\int_0^{u(t_i)}e^{H(t_i)}I_i(t){\rm d}t-\varepsilon\sum\limits_{i=0}^N\int_{s_i}^{t_{i+1}}e^{H(t)}G_i(t,u(t)){\rm d}t.
\end{eqnarray*}

\vskip2mm
\noindent
{\bf Lemma 2.1.} (\cite{Xia 2023,Heidarkhani}) {\it The embedding $E\hookrightarrow C^1([0,T])$ is compact and for all $u\in E$,  there hold
$$\|u\|_\infty\leq \|u\|\leq\dfrac{1}{\sqrt {e^{H_0}}}\|u\|_E,\qquad \|u'\|_\infty\leq \|u\|\leq\dfrac{1}{\sqrt {e^{H_0}}}\|u\|_E,$$
where $\|u\|_{\infty}=\max_{t\in [0,T]}|u(t)|$.}

\vskip2mm
\noindent
{\bf Lemma 2.2.} (\cite{Xia 2022} Lemma 2.2) {\it The embedding $E\hookrightarrow L^2([0,T])$ is compact and
$$\|u\|_2\leq\dfrac{T^2}{2\sqrt{2}}\|u\|\leq\dfrac{T^2}{2\sqrt {2e^{H_0}}}\|u\|_E,\qquad \|u^{\prime}\|_2\leq\dfrac{T}{\sqrt{2}}\|u\|\leq\dfrac{T}{\sqrt {2 e^{H_0}}}\|u\|_E,$$
for all $u\in E$, where $\|u\|_2=(\int_0^T|u(t)|^2{\rm d}t){^\frac{1}{2}}.$}

\vskip2mm
\noindent
{\bf Lemma 2.3.} (\cite{Kajikiya R 2013}) {\it Let $E$ be an infinite dimensional Banach space. For any $\varepsilon \in [0,1]$, ${J}_\varepsilon \in C(E,\mathbb{R})$. Suppose that ${J}_\varepsilon(u)$ has a continuous partial derivative ${J}_\varepsilon^\prime(u)$ with respect to $u$ and satisfies $(B_1)$-$(B_5)$ below.
\begin{itemize}
\item[$(B_1):$]$\inf\{{J}_\varepsilon(u):\varepsilon \in [0,1],u \in E\}> -\infty;$
\item[$(B_2):$]For $u \in E, \lvert {J}_\varepsilon(u)-{J}_0(u)\rvert\leq\psi(\varepsilon)$, where $\psi\in C([0,1], \mathbb{R})$ and $\psi(0)=0$, ${J}_0(u)=:{J}_\varepsilon(u)|_{\varepsilon=0}$;
\item[$(B_3):$]${J}_\varepsilon(u)$ satisfies the Palais-Smale condition uniformly on $\varepsilon$, i.e. if a sequence $\{(\varepsilon_k,u_k)\}\subset [0,1]\times E$ satisfies that
  $\sup_k\lvert {J}_{\varepsilon_k}(u_k)\rvert < \infty$ and ${J}_{\varepsilon_k}^\prime(u_k)$ converges to zero, then $\{(\varepsilon_k,u_k)\}$ has a convergent subsequence;
\item[$(B_4):$]${J}_0(u)={J}_0(-u)$ for $u\in E$ and ${J}_0(0)=0$;
\item[$(B_5):$]For any $u\in E\backslash \left\{0 \right\}$, there exists a unique $c(u)>0$ such that ${J}_0(cu)<0$ if $0<\lvert c\rvert <c(u)$ and ${J}_0(cu)\geq 0$ if $\lvert c \rvert \geq c(u)$.
\end{itemize}
Denote
$$ S^k:=\{t\in \mathbb{R}^{k+1}:|t|=1\},$$
$$
A_k:=\{\gamma \in C(S^k,E):\gamma \;\text{is odd}\},
$$
$$
d_k:=\inf_{\gamma\in A_k}\max_{t\in S^k}{{J}_0(\gamma(t))}.
$$
Let $k\in \mathbb{N}\backslash\left\{0 \right\}$ satisfying $d_k<d_{k+1}$. Then there exist two constants $\varepsilon_{k+1}, c_{k+1}$ such that $0< \varepsilon_{k+1}\leq 1$, $d_{k+1}\leq c_{k+1}<-\psi(\varepsilon)$ for $\varepsilon\in [0,\varepsilon_{k+1}]$, $J_\varepsilon(\cdot)$ has a critical value in the interval $[d_{k+1}-\psi(\varepsilon),c_{k+1}+\psi(\varepsilon)]$.}

\section{Proofs}
  \setcounter{equation}{0}
It is easy to see that such local assumptions $(F_0)$-$(F_5)$, $(I_1)$-$(I_3)$ and $(G)$ can not ensure that $J_\varepsilon(u)$ is well defined on $E$.
In order to overcome such difficulty, motivated by \cite{Kajikiya R 2013}, a cut-off technique is required. To be precise, we define $m(s)\in C^1(\mathbb{R},[0,1])$ as an even cut-off function satisfying $sm^{\prime}(s)\leq 0$ and strictly decreasing in $(\frac{\delta}{2},\delta)$, and satisfying
\begin{equation*}
m(s)=\left\{\begin{aligned}
1,&\quad \mbox{if} \; \lvert s\rvert \leq \frac{\delta}{2},\\
0,&\quad \mbox{if} \; \lvert s\rvert \geq\ \delta.
\end{aligned}\right.
\end{equation*}
Define $\tilde{F_i}:(s_i,t_{i+1}]\times \mathbb{R}\to \mathbb{R}$, $\tilde{G_i}:(s_i,t_{i+1}]\times \mathbb{R}\to \mathbb{R}$ and $\tilde{I_i}:\mathbb{R}\to \mathbb{R}$ as
\begin{eqnarray*}
&  &\tilde{F_i}(t,u):=m(u)F_i(t,u)+(1-m(u))a_*|u|^{\theta_1},\  i=0,1,...,N,\\
&  &\tilde{G_i}(t,u):=m(u)G_i(t,u),\ i=0,1,...,N,\\
&  &\tilde{I_i}(u):=m(u)I_i(u),\ i=1,2,...,N,
\end{eqnarray*}
Then, by direct calculation, we have
\begin{eqnarray*}
\tilde{f_i}(t,u)&:=& \frac{\partial}{\partial u}\tilde{F_i}(t,u)\\
&=&m'(u)F_i(t,u)+m(u)f_i(t,u)+a_*\theta_1|u|^{\theta_1-2}u-m'(u)a_*|u|^{\theta_1}-m(u)a_*\theta_1|u|^{\theta_1-2}u\\
&=&  m'(u)\left(F_i(t,u)-a_*|u|^{\theta_1}\right)+m(u)f_i(t,u)+a_*\theta_1(1-m(u))|u|^{\theta_1-2}u,\  i=0,1,...,N,
\end{eqnarray*}
and
$$
\tilde{g_i}(t,u):=\frac{\partial}{\partial u}m(u)G_i(t,u)=m'(u)G_i(t,u)+m(u)g_i(t,u),\  i=0,1,...,N.
$$

\par
Next, we give some properties for $\tilde{f_i}, \tilde{g_i}$ and $\tilde{I_i}$.

\vskip2mm
\noindent
{\bf Lemma 3.1.} {\it Assume that conditions $(F_0)$-$(F_5)$, $(I_1)$-$(I_3)$ and $(G)$ hold. Then $\tilde{f_i}(t,u)$ and $\tilde{g_i}(t,u)$ are continuous on $(s_i,t_{i+1}]\times \mathbb{R}$ and $\tilde{I_i}(u)$ is  continuous on $\mathbb{R}$, and the following conditions  hold:
\begin{itemize}
\item[$(F_1^{\prime})$]$\tilde{f_i}(t,u)=-\tilde{f_i}(t,-u), \mbox{ for all } (t,u)\in (s_i,t_{i+1}]\times \mathbb{R}, \ i=0,1,2,...,N$;
\item[$(F_2^{\prime})$]there hold\\
\hspace*{1.3cm} (i)\ $u\tilde{f}_i(t,u)-2\tilde{F}_i(t,u)<0$ if\  $(t,u)\in (s_i,t_{i+1}]\times \mathbb{R}\backslash\{0\}$, $i=0,1,...,N$,\\
\hspace*{1.3cm} (ii)\ $\tilde{g}_i(t,u)=\tilde{G}_i(t,u)=0$ if $|u|\geq \delta$ and $t\in (s_i,t_{i+1}]$, $i=0,1,...,N$,\\
\hspace*{1.3cm} (iii)\ $\tilde{F}_i(t,u)\geq0$ for all $(t,u)\in (s_i,t_{i+1}]\times \R$, $i=0,1,...,N;$
\item[$(F_3^{\prime})$]$\lim\limits_{|u| \to 0}\inf\limits_{t\in (s_i,t_{i+1}]}\dfrac {\tilde{F_{i}}(t,u)}{|u|^2}=+\infty,$   $i=0,1,2,...,N$;
\item[$(F_4^{\prime})$]there exists a constant $A_1>0$ such that
$$\lvert \tilde{f_i}(t,u)\rvert\le A_1\lvert u\rvert^{\theta_1-1}, \mbox{ for all }(t,u)\in (s_i,t_{i+1}]\times \mathbb{R},\ i=0,1,2,...,N;$$

\item[$(I_1^{\prime})$]$\tilde{I_i}(u)=-\tilde{I}_i(-u),$ for all $u\in \mathbb{R},\ i=1,2,...,N;$
\item[$(I_2^{\prime})$]there holds
$$\lvert \tilde{I_i}(u)\rvert\leq a_2(\lvert u\rvert^{\theta_2-1}+1),\; \text{for\;all}\;u\in \mathbb{R},\ i=1,2,...,N;$$

\item[$(I'_3)$] there holds\\
\hspace*{1.3cm} (i)\   $\hat{I}_i(u):=\int_0^{u} e^{H(t_i)} \tilde{I}_i(t){\rm d}t\geq0, \mbox{ for all } u\in \mathbb{R},\ i=1,2,...,N;$\\
\hspace*{1.3cm} (ii)\  $2\sum\limits_{i=1}^N\int_0^{u}e^{H(t_i)}\tilde{I}_i(t){\rm d}t\ge \sum\limits_{i=1}^Ne^{H(t_i)}\tilde{I}_i(u)u, \mbox{ for all } u\in \mathbb{R},\ i=1,2,...,N;$

\item[$(G')$] there exists a function $B\in C(\mathbb{R},\R^+)$ such that
    $$
    |\tilde{g}_i(t,u)|\le B(|u|), \mbox{ for all }t\in (s_i,t_{i+1}],\  i=0,1,2,...,N.
    $$

\end{itemize}}
\noindent
{\bf Proof.}\;By the definition of $\tilde{f_i}(t,u)$, $\tilde{g_i}(t,u)$ and $\tilde{I_i}(u)$, it is easy to obtain the continuity of  $\tilde{f_i}(t,u)$, $\tilde{g_i}(t,u)$ and $\tilde{I_i}(u)$.
 Next, we will prove that the conditions $(F_1^{\prime})$-$(F_4^{\prime})$, $(I_1^{\prime})$-$(I_3^{\prime})$, and $(G')$ hold in the subsequence.
 \par
{\bf $\bullet$}  If $|u|\leq\frac{\delta}{2}$, then $\tilde{f_i}(t,u)=f_i(t,u)$. Furthermore, in view of $(F_1)$ and the fact that $m(s)$ is an even function, we obtain that  $\tilde{f_i}(t,-u)=-\tilde{f_i}(t,u)$.
If $\frac{\delta}{2}<|u|<\delta$, by $(F_1)$ and the definition of  the function $m$, then
\begin{eqnarray*}
\tilde{f_i}(t,-u)&=& m'(-u)\left(F_i(t,-u)-a_*|-u|^{\theta_1}\right)+m(-u)f_i(t,-u)+a_*\theta_1(1-m(-u))|-u|^{\theta_1-2}(-u)\\
&=&-m'(u)\left(F_i(t,u)-a_*|u|^{\theta_1}\right)-m(u)f_i(t,u)-a_*\theta_1(1-m(u))|u|^{\theta_1-2}u\\
&=&-\tilde{f_i}(t,u).
\end{eqnarray*}
If $|u|\geq\delta$, then $\tilde{f_i}(t,u)=a_*\theta_1|u|^{\theta_1-2}u$ and then it is easy to obtain that $\tilde{f_i}(t,-u)=-\tilde{f_i}(t,u)$. Hence, $(F_1') $ holds.

\hspace*{0.1cm}{\bf $\bullet$} (i)\  If $0<|u|\le\frac{\delta}{2}$, then it follows from  $(F_2)$ that
\begin{eqnarray*}
u\tilde{f_i}(t,u)-2\tilde{F_i}(t,u)=uf_i(t,u)-2F_i(t,u)<0.
\end{eqnarray*}
If $\frac{\delta}{2}<|u|<\delta$,  then by the fact $s m^{\prime}(s)\leq 0$, $(F_2)$ and $(F_5)$, there holds
\begin{eqnarray*}
u\tilde{f_i}(t,u) -2\tilde{F_i}(t,u)
&=&  um'(u)\left(F_i(t,u)-a_*|u|^{\theta_1}\right)+m(u)(uf_i(t,u)-2F_i(t,u))\nonumber\\
&  &   +a_*(\theta_1-2)(1-m(u))|u|^{\theta_1}\nonumber\\
&<&  0.
\end{eqnarray*}
If $|u|\ge \delta$, then it follows from $1<\theta_1<2$ that
\begin{eqnarray*}
u\tilde{f_i}(t,u)-2\tilde{F_i}(t,u)=a_*(\theta_1-2)|u|^{\theta_1}<0.
\end{eqnarray*}

\par
(ii)\  If $|u|\geq\delta$,  then it follows from the fact $m(u)=0$ that
\begin{eqnarray*}
&  & \tilde{g}_i(t,u)=m'(u)G_i(t,u)+m(u)g_i(t,u)=0,\\
&  & \tilde{G}_i(t,u)=m(u)G_i(t,u)=0.
\end{eqnarray*}
\par
(iii)\  If $|u|\le\frac{\delta}{2}$ then it follows from ($F_2$) that $\tilde{F}_i(t,u)=F_i(t,u)\geq 0$. If $\frac{\delta}{2}<|u|<\delta$, it is easy to show that
\begin{eqnarray*}
\tilde{F}_i(t,u)=m(u)F_i(t,u)+(1-m(u))a_*|u|^{\theta_1}\geq0.
\end{eqnarray*}
If $|u|\geq\delta$, we have $\tilde{F}_i(t,u)=a_*|u|^{\theta_1}\geq0$. Thus $(F_2')(iii)$ holds.
\par

\hspace*{0.1cm}{\bf $\bullet$} If $|u| \to 0$, then we can make the restriction that $|u|<\frac{\delta}{2}$. Thus, $(F_3')$ can be obtained directly from $(F_3)$.

\hspace*{0.1cm}{\bf $\bullet$}  If $|u|\leq\frac{\delta}{2}$, then $(F_4)$ implies that $\lvert \tilde{f_i}(t,u)\rvert=\lvert f_i(t,u)\rvert\le a_1\lvert u\rvert^{\theta_1-1}$.
If $|u|\geq\delta$, it is obvious to see that $|\tilde{f}_i(t,u)|=a_*\theta_1|u|^{\theta_1-1}\le a_1\lvert u\rvert^{\theta_1-1}$.
When $\frac{\delta}{2}<|u|<\delta$,  by $(F_4)$, $(F_5)$ and the fact that $m(u)\in C^1(\mathbb{R},[0,1])$, we obtain
\begin{eqnarray}\label{eq1.2133}
|\tilde{f}_i(t,u)|&=&\Big|m'(u)\left(F_i(t,u)-a_*|u|^{\theta_1}\right)+m(u)f_i(t,u)+a_*\theta_1(1-m(u))|u|^{\theta_1-2}u\Big|\nonumber\\
&\leq&\dfrac{|um'(u)| \left(|F_i(t,u)|+a_*|u|^{\theta_1}\right)+|um(u)f_i(t,u)|+a_*\theta_1(1-m(u))|u|^{\theta_1}}{|u|}\nonumber\\
&\leq&\dfrac{|um'(u)|(\frac{a_1}{\theta_1}|u|^{\theta_1}+a_*|u|^{\theta_1})+|u m(u)f_i(t,u)|+a_*\theta_1(1-m(u))|u|^{\theta_1}}{|u|}\nonumber\\
&\leq&\dfrac{2A_0a_1|u|^{\theta_1}+|um(u)f_i(t,u)|+a_*\theta_1(1-m(u))|u|^{\theta_1}}{|u|}\nonumber\\
&\leq&2A_0a_1|u|^{\theta_1-1}+|f_i(t,u)|+a_*\theta_1|u|^{\theta_1-1}\nonumber\\
&\leq&2A_0a_1|u|^{\theta_1-1}+a_1|u|^{\theta_1-1}+a_*\theta_1|u|^{\theta_1-1}\nonumber\\
&\leq&2(A_0+1)a_1|u|^{\theta_1-1},\nonumber
\end{eqnarray}
where $A_0=\delta\cdot \max\limits_{s\in[\frac{\delta}{2}, \delta]}|m'(s)|$.  Let $A_1=(A_0+2)a_1$.  Thus $(F_4')$ holds.

\hspace*{0.1cm}{\bf $\bullet$} It is easy to obtain that $(I_1')$ holds from the definition of $\tilde{I}_i$ and the function $m$.

\hspace*{0.1cm}{\bf $\bullet$} Obviously, $(I_2') $ holds for all  $|u|\leq\frac{\delta}{2}$ and $|u|\geq\delta$.
If $\frac{\delta}{2}<|u|<\delta$, owing to the fact $m(s)\in C^1(\mathbb{R},[0,1])$, then
\begin{eqnarray*}
\lvert \tilde{I_i}(t,u)\rvert=\lvert m(u)I_i(t,u)\rvert\leq m(u)\cdot a_2(\lvert u\rvert^{\theta_2-1}+1)\leq a_2(\lvert u\rvert^{\theta_2-1}+1).
\end{eqnarray*}
Hence, $(I_2') $ holds for all $u\in \R$.

\hspace*{0.1cm}{\bf $\bullet$} (i)\  Since $m(s)$ is an even function, then  $(I_1)$ implies that $\hat{I}_i(u)=\int_0^{u}e^{H(t_i)}\tilde{I}_i(t){\rm d}t$ is even.
 Thus, without loss of generality, we can assume that $u\geq 0$. Thus, $(I_3)(i)$ shows that
\begin{eqnarray*}
\hat{I}_i(u)=\int_0^{u}e^{H(t_i)}\tilde{I}_i(t){\rm d}t=\int_0^{u}e^{H(t_i)}m(u)I_i(t)dt
=\begin{cases}
\int_0^{u}e^{H(t_i)}m(u)I_i(t)dt\ge 0,\ \ \mbox{ if } u\le \delta\\
\int_0^{\delta}e^{H(t_i)}m(u)I_i(t)dt\ge 0,\ \ \mbox{ if } u> \delta.
\end{cases}
\end{eqnarray*}
(ii)\ If $|u|\geq \delta$, then
\begin{eqnarray}\label{aa1}
\sum\limits_{i=1}^N e^{H(t_i)}\tilde{I}_i(u)u = \sum\limits_{i=1}^N e^{H(t_i)} m(u){I}_i(u)u =0.
\end{eqnarray}
Thus, by $(I_3)(i)$,  for all $u\geq\delta$,  we have
\begin{eqnarray}\label{aa2}
2\sum\limits_{i=1}^N\int_0^{u}e^{H(t_i)}\tilde{I}_i(t){\rm d}t
& =&  2\sum\limits_{i=1}^N\int_0^{\frac{\delta}{2}}e^{H(t_i)}\tilde{I}_i(t){\rm d}t+2\sum\limits_{i=1}^N\int_{\frac{\delta}{2}}^{ \delta}e^{H(t_i)}\tilde{I}_i(t){\rm d}t +2\sum\limits_{i=1}^N\int_{ \delta}^{u}e^{H(t_i)}\tilde{I}_i(t){\rm d}t\nonumber\\
&=&  2\sum\limits_{i=1}^N\int_0^{\frac{\delta}{2}}e^{H(t_i)}{I}_i(t){\rm d}t+2\sum\limits_{i=1}^N\int_{\frac{\delta}{2}}^{ \delta}e^{H(t_i)}m(t){I}_i(t){\rm d}t \nonumber\\
&\geq& 2\sum\limits_{i=1}^N\int_0^{\frac{\delta}{2}}e^{H(t_i)}{I}_i(t){\rm d}t\nonumber\\
& \ge & 0,
\end{eqnarray}
which, together with (\ref{aa1}) and (\ref{aa2}), implies that $2\sum\limits_{i=1}^N\int_0^{u}e^{H(t_i)}\tilde{I}_i(t){\rm d}t\ge \sum\limits_{i=1}^Ne^{H(t_i)}\tilde{I}_i(u)u$ for all $|u|\geq \delta$. If $|u|< \delta$, by $(I_3)(ii)$ and the fact that $m(s)\in [0,1]$ and is strictly decreasing in $(\frac{\delta}{2},\delta)$, we have
\begin{eqnarray*}
2\sum\limits_{i=1}^N\int_0^{u}e^{H(t_i)}\tilde{I}_i(t){\rm d}t
 =  2\sum\limits_{i=1}^N\int_0^{u}e^{H(t_i)}m(t){I}_i(t){\rm d}t
\ge  2\sum\limits_{i=1}^N\int_0^{u}e^{H(t_i)}m(u){I}_i(t){\rm d}t
\ge  \sum\limits_{i=1}^N e^{H(t_i)} m(u){I}_i(u)u.
\end{eqnarray*}

\hspace*{0.1cm}{\bf $\bullet$} If $|u|\leq\frac{\delta}{2}$, we have $|\tilde{g_i}(t,u)|=|m'(u)G_i(t,u)+m(u)g_i(t,u)|=|g_i(t,u)|\le b(|u|)$. If $|u|\geq\delta$, $\tilde{g_i}(t,u)=0\le b(|u|)$. If $\frac{\delta}{2}<|u|<\delta$, then
\begin{eqnarray*}
|\tilde{g_i}(t,u)|&=&|m'(u)G_i(t,u)+m(u)g_i(t,u)|\\
&\leq& |m'(u)G_i(t,u)|+|m(u)g_i(t,u)|\\
&\leq& |m'(u)\int_0^u g_i(t,s)ds|+|g_i(t,u)|\\
&\leq& |m'(u)|\int_0^{|u|} |g_i(t,s)|ds+|g_i(t,u)|\\
&\leq& K|u|\cdot \max_{s\in [0,\delta]}b(s)+b(|u|)\\
&\leq& K\delta\cdot \max_{s\in [0,\delta]}b(s)+b(|u|)
\end{eqnarray*}
where $K=\max\limits_{u\in[\frac{\delta}{2}, \delta]}|m'(u)|$. Let $B(|u|)=K\delta\cdot \max_{s\in [0,\delta]}b(s)+b(|u|)$. Thus $(G')$ holds.
The proof is complete.
\qed

\vskip2mm
 \noindent
{\bf Remark 3.1.} By  $(F_4')$ and the fact
\begin{eqnarray*}
\tilde{F}_i(t,u)=\int_0^u \tilde{f}_i(t,s)ds,\ \forall (t,u)\in (s_i,t_{i+1}]\times \mathbb{R},\ i=0,1,...,N,
\end{eqnarray*}
it is easy to see that
\begin{eqnarray*}
|\tilde{F}_i(t,u)|\leq \frac{A_1}{\theta_1}|u|^{\theta_1},\ \forall (t,u)\in (s_i,t_{i+1}]\times \mathbb{R},\ i=0,1,...,N.
\end{eqnarray*}

\vskip2mm
 \noindent
{\bf Remark 3.2.} The condition $(F_2)'(i)$ implies the fact that $u\tilde{f}_i(t,u)-2\tilde{F}_i(t,u)=0$ if and only if $u=0$, $i=0,1,...,N$.
Indeed, if $u=0$, it is obvious that  $u\tilde{f}_i(t,u)-2\tilde{F}_i(t,u)=0$. Conversely, assume that  $u\tilde{f}_i(t,u)-2\tilde{F}_i(t,u)=0$ but $u\not=0$. This would contradicts $(F_2')(i)$.

\vskip2mm
\noindent
\par
Next, we consider the following modified model:
\begin{eqnarray}\label{eq1.23}
\begin{cases}
 u^{(4)}(t)+2h(t)u^{\prime \prime \prime}(t)+\left(h^{2}(t)+h^{\prime}(t)+\beta e^{-\int_0^t h(\tau)}\right)u^{\prime \prime}(t)=\tilde{f}_i(t,u(t))+\varepsilon \tilde{g}_i(t,u(t)),\; t\in (s_i,t_{i+1}],\\
 i=0,1,2,...,N,\\
 u^{(4)}(t)+2h(t)u^{\prime \prime \prime}(t)+\left(h^{2}(t)+h^{\prime}(t)+\beta e^{-\int_0^t h(\tau)}\right)u^{\prime \prime}(t)=0,\; t\in (t_i,s_i],\; i=1,2,...,N,\\
 u^{\prime \prime}(t_i^+)=u^{\prime \prime}(t_i^-),\ s^{\prime \prime}(t_i^+)=s^{\prime \prime}(t_i^-),\ u^{\prime \prime \prime}(s_i^+)=u^{\prime \prime \prime}(s_i^-),\; i=1,2,...,N,\\
 \Delta u^{\prime \prime \prime}(t_i)=\tilde{I}_i(u(t_i)),\; i=1,2,...,N,\\
 u(0)=u^{\prime}(0)=u(T)=u^{\prime}(T)=0,
 \end{cases}
\end{eqnarray}
and define the corresponding energy functional $\tilde{J}_\varepsilon(u)$ on $E$ by
\begin{eqnarray*}
\tilde{J}_\varepsilon(u)&:=&\frac{1}{2}\|u\|_E^2-\frac{\beta}{2}\int_0^T|u^{\prime}(t)|^2{\rm d}t-\sum\limits_{i=0}^N\int_{s_i}^{t_{i+1}}e^{H(t)}\tilde{F}_i(t,u(t)){\rm d}t\notag\\
&  & -\sum\limits_{i=1}^N\int_0^{u(t_i)}e^{H(t_i)}\tilde{I}_i(t){\rm d}t-\varepsilon\sum\limits_{i=0}^N\int_{s_i}^{t_{i+1}}e^{H(t)}\tilde{G}_i(t,u(t)){\rm d}t.
\end{eqnarray*}
By Lemma 2.1, Lemma 2.2 and the continuity of $\tilde{f_i}(t,u)$, $\tilde{g_i}(t,u)$ and $\tilde{I_i}(u)$, a standard argument implies that $\tilde{J}_\varepsilon$ is well-defined and $\tilde{J}_\varepsilon$ is continuously differentiable on $E$. Moreover,
\begin{eqnarray}\label{eq7}
\langle\tilde{J}_\varepsilon^{\prime}(u),v\rangle&=&\int_0^Te^{H(t)}u^{\prime\prime}(t)v^{\prime\prime}(t){\rm d}t-\beta\int_0^Tu^{\prime}(t)v^{\prime}(t){\rm d}t-\sum\limits_{i=0}^N\int_{s_i}^{t_{i+1}}e^{H(t)}\tilde{f}_i(t,u(t))v(t){\rm d}t\notag\\
&  & -\sum\limits_{i=1}^Ne^{H(t_i)}\tilde{I}_i(u(t_i))v(t_i)-\varepsilon\sum\limits_{i=0}^N\int_{s_i}^{t_{i+1}}e^{H(t)}\tilde{g}_i(t,u(t))v(t){\rm d}t
\end{eqnarray}
for any $u,v\in E$.
Thus, the critical points of $\tilde{J}_\varepsilon$ in $E$ are the weak solutions of system (\ref{eq1.23}) (see \cite{Xia 2023}).

\vskip2mm
 \noindent
{\bf Lemma 3.2.} {\it Assume that $(F'_4)$ and $(I'_2)$ hold. Then for all $u\in E$,
\begin{eqnarray*}
\sum\limits_{i=0}^N\int_{s_i}^{t_{i+1}}e^{H(t)}\tilde{F}_i(t,u(t)){\rm d}t\leq M_1A_1Te_*\|u\|_E^{\theta_1},
\end{eqnarray*}
and
\begin{eqnarray}\label{eq3.1}
\sum\limits_{i=1}^N\int_0^{u(t_i)}e^{H(t_i)}\tilde{I}_i(t){\rm d}t \leq NM_1a_2\cdot\frac{\theta_2 +1}{\theta_2} \cdot e_*(\|u\|_E^{\theta_2}+1),
\end{eqnarray}}
where $e_*=\max\left\{{{e^{-{H_{0}\theta_1 \over 2}}}},{e^{-{{H_{0}\theta_2\over 2}}}}\right\}$ and $ M_1=\max_{t\in [0,T]}{e^{H(t)}}$.
\vskip2mm
\noindent
{\bf Proof.}\;By Remark 3.1, Young's inequality and Lemma 2.1, we have
\begin{eqnarray*}
\sum\limits_{i=0}^N\int_{s_i}^{t_{i+1}}e^{H(t)}\tilde{F}_i(t,u(t)){\rm d}t&\leq&\sum\limits_{i=0}^N\int_{s_i}^{t_{i+1}}e^{H(t)}\frac{A_1}{\theta_1}|u(t)|^{\theta_1}{\rm d}t\\
&\leq& M_1A_1T\|u\|_{\infty}^{\theta_1}\\
&\leq& M_1A_1Te^{-\frac{H_0\theta_1}{2}}\|u\|_E^{\theta_1}\\
&\leq& M_1A_1Te_*\|u\|_E^{\theta_1},
\end{eqnarray*}
Similarly, by $(I'_2)$, it is easy to prove that (\ref{eq3.1}) is also valid, which can also be seen in \cite{Xia 2023}.
\vskip2mm
\par
Next, we prove that the functional $\tilde{J}_\varepsilon$ satisfy all the conditions in Lemma 2.3.
\vskip2mm
\noindent
{\bf Lemma 3.3.}  {\it Assume that $(F_2^\prime)$, $(F_4^\prime)$, $(I_2^\prime)$ and $(H)$ hold. Then the functional $\tilde{J}_\varepsilon$ is coercive and bounded below on $E$ for all $\varepsilon\in[0,1]$.}
\vskip2mm
\noindent
{\bf Proof.}\; For any  $u\in E$ and $\varepsilon\in[0,1]$, it follows from $(F'_2)$ that
\begin{eqnarray}
\Bigg |\varepsilon\int_{s_i}^{t_i+1}\tilde{G}_i(t,u(t))dt\Bigg |&\leq& |\varepsilon|\int_{\{t\in[s_i,t_{i+1}]:|u(t)|<\delta\}}|\tilde{G}_i(t,u(t))|{\rm d}t+|\varepsilon|\int_{\{t\in[s_i,t_{i+1}]:|u(t)|\ge\delta\}}|\tilde{G}_i(t,u(t))|{\rm d}t\nonumber\\
\label{eqq3.33}&= &|\varepsilon|\int_{\{t\in[s_i,t_{i+1}]:|u(t)|<\delta\}}|\tilde{G}_i(t,u(t))|{\rm d}t\\
\label{eq3.33}&\leq&\max_{(t,s)\in [0,T]\times[-\delta,\delta]}|\tilde{G}_i(t,s)|T.
\end{eqnarray}
Thus, for all $u\in E$ and $\beta>0$, it follows from Lemma 2.2, Lemma 3.2   and (\ref{eq3.33}) that
\begin{eqnarray}\label{eq3.34}
\tilde{J}_\varepsilon(u)&=&\frac{1}{2}\|u\|_E^2-\frac{\beta}{2}\int_0^T|u^{\prime}(t)|^2{\rm d}t-\sum\limits_{i=0}^N\int_{s_i}^{t_{i+1}}e^{H(t)}\tilde{F}_i(t,u(t)){\rm d}t\nonumber\\
& &-\sum\limits_{i=1}^N\int_0^{u(t_i)}e^{H(t_i)}\tilde{I}_i(t){\rm d}t-\varepsilon\sum\limits_{i=0}^N\int_{s_i}^{t_{i+1}}e^{H(t)}\tilde{G}_i(t,u(t)){\rm d}t\nonumber\\
&\geq& \left(\frac{1}{2}-\frac{\beta T^2}{4e^{H_0}}\right)\|u\|_E^2-M_1A_1Te_*\|u\|_E^{\theta_1}-NM_1a_2\cdot\frac{\theta_2 +1}{\theta_2} \cdot e_*\|u\|_E^{\theta_2}\nonumber\\
&  & -NM_1a_2\cdot \dfrac{\theta_2+1}{\theta_2}\cdot e_*-\max_{(t,s)\in [0,T]\times[-\delta,\delta]}|\tilde{G}_i(t,s)|T(N+1).
\end{eqnarray}
Similarly, for $\beta\le 0$, together with (\ref{eq3.34}), it is easy to see that
\begin{eqnarray}\label{eq3.35}
\tilde{J}_\varepsilon(u)&\geq&  \frac{1}{2}\|u\|_E^2-M_1A_1Te_*\|u\|_E^{\theta_1}-NM_1a_2\cdot\frac{\theta_2 +1}{\theta_2} \cdot e_*\|u\|_E^{\theta_2}\nonumber\\
&&   -NM_1a_2\cdot \dfrac{\theta_2+1}{\theta_2}\cdot e_*-\max_{(t,s)\in [0,T]\times[-\delta,\delta]}|\tilde{G}_i(t,s)|T(N+1).
\end{eqnarray}
Combining $(H)$, (\ref{eq3.34}), (\ref{eq3.35}), and the fact that $1<\theta_1<2$ and $1<\theta_2<2$, we deduce that $\tilde{J}_\varepsilon$ is coercive and then bounded below on $E$. \qed

\vskip2mm
\noindent
{\bf Lemma 3.4.} {\it Suppose that ($F_4^{\prime}$), ($I_2^{\prime}$), $(G')$ and $(H)$ hold.  Then $\tilde{J}_\varepsilon$ satisfies condition $(B_3)$ in Lemma 2.3.}
\vskip2mm
\noindent
{\bf Proof.}\;Let $\{(\varepsilon_n, u_n)\}\subset [0,1]\times E$ be any sequence such that
\begin{eqnarray}\label{a1}
\tilde{J}_{\varepsilon_n}(u_n)<\infty \;\text{and}\; \tilde{J}_{\varepsilon_n}^{\prime}(u_n)\rightarrow 0\;\text{as}\;n\to\infty.
\end{eqnarray}
Then, from (\ref{a1}) and Lemma 3.3, it follows that the sequence $\{(\varepsilon_n, u_n)\}$ is  bounded in $[0,1]\times E$ for all $\beta<\frac{2e^{H_0}}{T^2}$. Hence,
there exists a positive constant $M_0>0$ such that $\|u_n\|_E\le M_0$. Therefore, together with Lemma 2.1,  there exists a subsequence of $\{\varepsilon_n\}$, still denoted by $\{\varepsilon_n\}$, such that ${\varepsilon_n}$ converges to $\varepsilon$ in $[0,1]$, and a subsequence of $\{u_n\}$, still denoted by $\{u_n\}$, such that
\begin{eqnarray}\label{eq112}
u_n\rightharpoonup u\;\text{in}\;E,\qquad u_n\rightarrow u\;\text{in}\;C^1([0,T]).
\end{eqnarray}
Then it follows from $(G')$, (\ref{eq112}) and Lemma 2.1 that
\begin{eqnarray}\label{eqqq9}
&     & \sum_{i=0}^N\int_{s_i}^{t_{i+1}}e^{H(t)}|\tilde{g}_i(t,u_n(t))-\tilde{g}_i(t,u(t))||u_n(t)-u(t)|dt\nonumber\\
& \le & \sum_{i=0}^N\int_{s_i}^{t_{i+1}}e^{H(t)}(B(|u_n(t)|)+B(|u(t)|))||u_n(t)-u(t)|dt\nonumber\\
& \le & \|u_n-u\|_\infty M_1 T(N+1)\left(\max_{|s|\le\frac{M_0}{\sqrt{e^{H_0}}}} B(s)+\max_{|s|\le\|u\|_\infty} B(s)\right)\nonumber\\
&\to  & 0, \ \ \mbox{as }n\to \infty.
\end{eqnarray}
Similarly, by $(F_4')$ and $(I_2')$, we also have
\begin{eqnarray}\label{eq9}
\sum_{i=0}^N\int_{s_i}^{t_{i+1}}e^{H(t)}|\tilde{f}_i(t,u_n(t))-\tilde{f}_i(t,u(t))||u_n(t)-u(t)|dt=o(1)\; \text{as} \; n\rightarrow \infty
\end{eqnarray}
and
\begin{equation}\label{eq10}
\sum_{i=1}^Ne^{H(t_i)}|\tilde{I}_i(t,u_n(t_i))-\tilde{I}_i(t,u(t_i))||u_n(t_i)-u(t_i)|=o(1)\; \text{as} \; n\rightarrow \infty.
\end{equation}
Note that if  $\beta>0$, by (\ref{eq7}), then
\begin{eqnarray}\label{eq113}
&&\langle \tilde{J}'_{\varepsilon_n}(u_n)-\tilde{J}'_{\varepsilon}(u), u_n-u \rangle\nonumber\\
&=&\langle \tilde{J}'_{\varepsilon_n}(u_n),  u_n-u \rangle-\langle\tilde{J}'_{\varepsilon}(u), u_n-u \rangle\nonumber\\
&=&\int_0^Te^{H(t)}u_n^{\prime\prime}(u_n^{\prime\prime}-u^{\prime\prime}){\rm d}t-\beta\int_0^Tu_n^{\prime}(u_n^{\prime}-u^{\prime}){\rm d}t-\sum\limits_{i=0}^N\int_{s_i}^{t_{i+1}}e^{H(t)}\tilde{f}_i(t,u_n)(u_n-u){\rm d}t\nonumber\\
&    &-\sum\limits_{i=1}^Ne^{H(t_i)}\tilde{I}_i(u_n(t_i))(u_n(t_i)-u(t_i))-\varepsilon_n\sum\limits_{i=0}^N\int_{s_i}^{t_{i+1}}e^{H(t)}\tilde{g}_i(t,u_n)(u_n-u){\rm d}t\nonumber\\
&    &-\Bigg[\int_0^Te^{H(t)}u^{\prime\prime}(u_n^{\prime\prime}-u^{\prime\prime}){\rm d}t-\beta\int_0^Tu^{\prime}(u_n^{\prime}-u^{\prime}){\rm d}t{\rm d}t-\sum\limits_{i=0}^N\int_{s_i}^{t_{i+1}}e^{H(t)}\tilde{f}_i(t,u)(u_n-u){\rm d}t\nonumber\\
&    &-\sum\limits_{i=1}^Ne^{H(t_i)}\tilde{I}_i(u(t_i))(u_n(t_i)-u(t_i))-\varepsilon\sum\limits_{i=0}^N\int_{s_i}^{t_{i+1}}e^{H(t)}\tilde{g}_i(t,u)(u_n-u){\rm d}t\Bigg]\nonumber\\
&\geq&\|u_n-u\|_E^2-\beta\int_0^T(u_n^{\prime}-u^{\prime})^2{\rm d}t-\sum\limits_{i=0}^N\int_{s_i}^{t_{i+1}}e^{H(t)}|\tilde{f}_i(t,u_n)-\tilde{f}_i(t,u)||u_n-u|{\rm d}t\nonumber\\
&    &-\sum\limits_{i=1}^Ne^{H(t_i)}|\tilde{I}_i(u_n(t_i))-\tilde{I}_i(u(t_i))||u_n(t_i)-u(t_i)|\nonumber\\
&    & -|\varepsilon_n-\varepsilon|\sum\limits_{i=0}^N\int_{s_i}^{t_{i+1}}e^{H(t)}|\tilde{g}_i(t,u_n)-\tilde{g}_i(t,u)||u_n-u|{\rm d}t\nonumber\\
&\geq& \|u_n-u\|_E^2-\frac{\beta T^2}{2e^{H_0}}\|u_n-u\|_E^2-\sum\limits_{i=0}^N\int_{s_i}^{t_{i+1}}e^{H(t)}|\tilde{f}_i(t,u_n)-\tilde{f}_i(t,u)||u_n-u|{\rm d}t\nonumber\\
&    & -\sum\limits_{i=1}^Ne^{H(t_i)}|\tilde{I}_i(u_n(t_i))-\tilde{I}_i(u(t_i))||u_n(t_i)-u(t_i)|\nonumber\\
&    & -|\varepsilon_n-\varepsilon|\sum\limits_{i=0}^N\int_{s_i}^{t_{i+1}}e^{H(t)}|\tilde{g}_i(t,u_n)-\tilde{g}_i(t,u)||u_n-u|{\rm d}t,
\end{eqnarray}
If $\beta \leq 0$, then we have
\begin{eqnarray}\label{eq113a}
&      & \langle \tilde{J}'_{\varepsilon_n}(u_n)-\tilde{J}'_{\varepsilon}(u), u_n-u \rangle\nonumber\\
& \geq & \|u_n-u\|_E^2-\sum\limits_{i=0}^N\int_{s_i}^{t_{i+1}}e^{H(t)}|\tilde{f}_i(t,u_n)-\tilde{f}_i(t,u)||u_n-u|{\rm d}t-\sum\limits_{i=1}^Ne^{H(t_i)}|\tilde{I}_i(u_n(t_i)) -\tilde{I}_i(u(t_i))||u_n(t_i)-u(t_i)|\nonumber\\
&      & -|\varepsilon_n-\varepsilon|\sum\limits_{i=0}^N\int_{s_i}^{t_{i+1}}e^{H(t)}|\tilde{g}_i(t,u_n)-\tilde{g}_i(t,u)||u_n-u|{\rm d}t.
\end{eqnarray}
Thus together with  (\ref{a1}), (\ref{eq112}), (\ref{eqqq9}), (\ref{eq9}), (\ref{eq10}), (\ref{eq113}) and  (\ref{eq113a}) imply that
\begin{eqnarray*}
o(1)&=&\Big| \langle \tilde{J}'_{\varepsilon_n}(u_n)-\tilde{J}'_{\varepsilon}(u), u_n-u \rangle \Big|\\
&\geq&\begin{cases} (1-\frac{\beta T^2}{2e^{H_0}})\|u_n-u\|_E^2+o(1), \ \mbox{if } \beta>0,\\
\|u_n-u\|_E^2+o(1),\ \mbox{if }\beta\le 0,
\end{cases}\ \ \text{as} \; n\rightarrow \infty,
\end{eqnarray*}
which shows that $\|u_n-u\|_E^2\rightarrow 0$ as $n\rightarrow \infty$ by $(H)$.
\qed

\vskip2mm
\noindent
{\bf Lemma 3.5.} {\it For any given $u\in E\backslash \{0\}$, $\tilde{J}_0(u)$ satisfies condition $(B_5)$ in Lemma 2.3.}
\vskip2mm
\noindent
{\bf Proof.}  For any given $u\in E\backslash \{0\}$, by $(F_1')$ and $(I_1')$, we may assume $c>0$ and let
$$
P(c)=\frac{1}{2}\|u\|_E^2-\frac{\beta}{2}\int_0^T|u^{\prime}(t)|^2{\rm d}t-c^{-2}\sum\limits_{i=0}^N\int_{s_i}^{t_{i+1}}e^{H(t)}\tilde{F}_i(t,cu(t)){\rm d}t,
$$
and
$$
Q(c)=-c^{-2}\sum\limits_{i=1}^N\int_0^{cu(t_i)}e^{H(t_i)}\tilde{I}_i(t){\rm d}t=-c^{-2}\sum\limits_{i=1}^N\hat{I}_i(cu(t_i)).
$$
Then
$$
\tilde{J}_0(cu)=c^2(P(c)+Q(c)).
$$
Define $\Upsilon_u(c):=P(c)+Q(c)$. Next, we will  prove the existence and uniqueness of $c(u)$, respectively.

{\bf Existence.} Since $u(t)\not\equiv 0$ for all $t\in [0,T]$, by the continuity of $u$, there exists at least an interval $(a_u,b_u)\subset (s_{i_0}, t_{i_0+1}]$ for some
 $i_0\in \{1,...,N\}$ such that $u(t)\not= 0$ for all $t\in (a_u,b_u)$. Thus, we can take $\eta>0$ small enough such that
$$
 K_\eta:=\{t\in (s_{i_0},t_{i_0+1}]: \eta<|u(t)|<\frac{1}{\eta}\}\  \text{and}\; \mu(K_\eta)>0,
$$
where $\mu$ denotes the Lebesgue measure of $[0,T]$. Next, we shall discuss the cases that $\beta> 0$ and $\beta\leq0$, respectively.

{\rm Case 1.}  For $\beta>0$, due to $(F_2')(iii)$ and $(I'_3)$, we can estimate the function $\Upsilon_u(c)$ as follows: for $0<c<1$,
\begin{eqnarray}\label{eq114}
\Upsilon_u(c)&=&\frac{1}{2}\|u\|_E^2-\frac{\beta}{2}\int_0^T|u^{\prime}(t)|^2{\rm d}t-c^{-2}\sum\limits_{i=0}^N\int_{s_i}^{t_{i+1}}e^{H(t)}\tilde{F}_i(t,cu(t)){\rm d}t-c^{-2}\sum\limits_{i=1}^N\int_0^{cu(t_i)}e^{H(t_i)}\tilde{I}_i(t){\rm d}t\nonumber\\
&\leq& \frac{1}{2}\|u\|_E^2-c^{-2}\sum\limits_{i=0}^N\int_{s_i}^{t_{i+1}}e^{H(t)}\tilde{F}_i(t,cu(t)){\rm d}t\nonumber\\
&\leq& \frac{1}{2}\|u\|_E^2-c^{-2}\int_{s_{i_0}}^{t_{i_0+1}}e^{H(t)}\tilde{F}_{i_0}(t,cu(t)){\rm d}t\nonumber\\
&\leq&\frac{1}{2}\|u\|_E^2-\int_{K_\eta}\frac{e^{H(t)}\tilde{F}_{i_0}(t,cu(t))|u(t)|^2}{|cu(t)|^2}{\rm d}t\nonumber\\
&\leq& \frac{1}{2}\|u\|_E^2-\inf_{t\in K_\eta}\frac{\tilde{F}_{i_0}(t,cu(t))}{|cu(t)|^2}\int_{K_\eta}e^{H(t)}|u(t)|^2dt\nonumber\\
&\leq& \frac{1}{2}\|u\|_E^2-M_1\eta^2\mu(K_{\eta})\inf_{t\in K_\eta}\frac{\tilde{F}_{i_0}(t,cu(t))}{|cu(t)|^2}.
\end{eqnarray}
By condition $(F_3^{\prime})$, (\ref{eq114}) implies that $\lim\limits_{c\rightarrow 0^+}\Upsilon_u(c)=-\infty$. Hence, $\Upsilon_u(c)<0$ for $c>0$ small enough.
On the other hand, for $c>0$ large enough, by using Lemma 2.2, Lemma 3.2 and $(H)$, we have
\begin{eqnarray*}
\Upsilon_u(c)&=&\frac{1}{2}\|u\|_E^2-\frac{\beta}{2}\int_0^T|u^{\prime}(t)|^2{\rm d}t-c^{-2}\sum\limits_{i=0}^N\int_{s_i}^{t_{i+1}}e^{H(t)}\tilde{F}_i(t,cu(t)){\rm d}t-c^{-2}\sum\limits_{i=1}^N\int_0^{cu(t_i)}e^{H(t_i)}\tilde{I}_i(t){\rm d}t\nonumber\\
&\geq&\frac{1}{2}\|u\|_E^2-\frac{\beta}{2}\|u'\|_2^2-c^{-2}M_1A_1Te_*\|cu\|_E^{\theta_1}-c^{-2}NM_1a_2\cdot\frac{\theta_2 +1}{\theta_2} \cdot e_*(\|cu\|_E^{\theta_2}+1)\nonumber\\
&=&\left(\frac{1}{2}-\frac{\beta T^2}{4e^{H_0}}\right)\|u\|_E^2-c^{{\theta_1}-2}M_1A_1Te_*\|u\|_E^{\theta_1}-c^{\theta_2-2}NM_1a_2\cdot\frac{\theta_2 +1}{\theta_2} \cdot e_*(\|u\|_E^{\theta_2}+1)\nonumber\\
&>& 0,
\end{eqnarray*}
where $1<\theta_1<2$ and $1<\theta_2<2$.

{\rm Case 2.}  For $\beta\leq0$, by  using the same proofs as (\ref{eq114}), for $0<c<1$, we have
\begin{eqnarray*}
\Upsilon_u(c)&\leq&\left(\frac{1}{2}-\frac{\beta T^2}{4e^{H_0}}\right) \|u\|_E^2-M_1\eta^2\mu(K_{\eta})\inf_{t\in K_\eta}\frac{\tilde{F}_{i_0}(t,cu(t))}{|cu(t)|^2},
\end{eqnarray*}
which, together with $(F_3)$, implies that $\lim\limits_{c\rightarrow 0^+}\Upsilon_u(c)=-\infty$. Hence, $\Upsilon_u(c)<0$ for $c>0$ small enough. On the other hand, for $c>0$ large enough, we have
\begin{eqnarray*}
\Upsilon_u(c)\geq\frac{1}{2}\|u\|_E^2-c^{\theta_1-2}M_1A_1Te_*\|u\|_E^{\theta_1}-c^{\theta_2-2}NM_1a_2\cdot\frac{\theta_2 +1}{\theta_2} \cdot e_*(\|u\|_E^{\theta_2}+1)>0.
\end{eqnarray*}

{\bf Uniqueness.} By  condition $(F_2^{\prime})$, we have
$$
P^{\prime}(c)=-c^{-3}\sum\limits_{i=0}^N\int_{s_i}^{t_{i+1}}e^{H(t)}\big(cu\tilde{f}_i(t,cu)-2\tilde{F}_i(t,cu)\big){\rm d}t> 0\;\text{for all}\; c>0,
$$
and by condition (\ref{eq3.1}) and $(I_3^{\prime})$, we have
\begin{eqnarray*}
Q^{\prime}(c)&=&2c^{-3}\sum\limits_{i=1}^N\hat{I}_i(cu(t_i))-c^{-3}\sum\limits_{i=1}^N  e^{H(t_i)}  cu(t_i)\tilde{I}_i(cu(t_i))\\
&=&c^{-3}\bigg(2\sum\limits_{i=1}^N\hat{I}_i(cu(t_i))-\sum\limits_{i=1}^N  e^{H(t_i)}  cu(t_i)\tilde{I}_i(cu(t_i))\bigg)\\
&\ge & 0,
\end{eqnarray*}
for all $c>0$. Thus, $\Upsilon_u^{\prime}(c)>0$ on  $(0,+\infty)$, that is, $\Upsilon_u$ is strictly increasing on $ (0,+\infty)$.
\par
To sum up, for any given $u\in E\backslash \{0\}$ and $\beta<\frac{2e^{H_0}}{T^2}$, $\Upsilon_u(c)$ has a unique zero $c(u)$ such that
$$
\Upsilon_u(c)<0\; \text{for}\; 0<c<c(u),\quad \mbox{ and } \Upsilon_u(c)\geq 0  \quad  \text{for all } c\geq c(u).
$$
The proof is complete.
\qed

 \vskip2mm
 \noindent
{\bf Lemma 3.6.}  {\it For any $b>0$, there exists a $\sigma(b)>0$ such that if  $\tilde{J}_\varepsilon^{\prime}(u)=0$ and $|\tilde{J}_\varepsilon(u)|\leq\sigma(b)$, where $|\varepsilon|\leq\sigma(b)$, then $\|u\|_E\leq b$.}
\vskip2mm
\noindent
{\bf Proof.}  Assume by contradiction that there exist two sequences $\{u_n\}\subset E$ and $\{\varepsilon_n\}$ with $\varepsilon_n\to 0$  as $n\to \infty$, satisfying $\tilde{J}_{\varepsilon_n}(u_n)\to 0$, $\tilde{J}_{\varepsilon_n}^{\prime}(u_n)=0$ and $\|u_n\|_E\geq b_0>0$ for all $n$, where $b_0$ is a constant independent of $n$. Then by Lemma 3.4, there exists a subsequence of $\{u_n\}$ that converges to $u_0$ in $E$ and then
\begin{eqnarray*}
0=\langle\tilde{J}_0'(u_0),u_0\rangle& = &\int_0^Te^{H(t)}|u_0^{\prime\prime}(t)|^2{\rm d}t-\beta\int_0^T|u_0^{\prime}(t)|^2{\rm d}t-\sum\limits_{i=0}^N\int_{s_i}^{t_{i+1}}e^{H(t)}\tilde{f}_i(t,u_0(t))u_0(t){\rm d}t\\
&  &-\sum\limits_{i=1}^Ne^{H(t_i)}\tilde{I}_i(u_0(t_i))u_0(t_i)
\end{eqnarray*}
and
\begin{eqnarray*}
0=\tilde{J}_0(u_0)=\frac{1}{2}\|u_0\|_E^2-\frac{\beta}{2}\int_0^T|u_0^{\prime}(t)|^2{\rm d}t-\sum\limits_{i=0}^N\int_{s_i}^{t_{i+1}}e^{H(t)}\tilde{F}_i(t,u_0(t)){\rm d}t -\sum\limits_{i=1}^N\int_0^{u_0(t_i)}e^{H(t_i)}\tilde{I}_i(t){\rm d}t.
\end{eqnarray*}
From the above two equations, it holds
\begin{eqnarray}\label{eq1.156}
0&=&2\tilde{J}_0(u_0)-\langle\tilde{J}_0'(u_0),u_0\rangle\nonumber\\
&=& \sum\limits_{i=0}^N\int_{s_i}^{t_{i+1}}e^{H(t)}(\tilde{f}_i(t,u_0(t))u_0(t)-2\tilde{F}_i(t,u_0(t))){\rm d}t \nonumber\\
&  & +\sum\limits_{i=1}^Ne^{H(t_i)}\tilde{I}_i(u_0(t_i))u_0(t_i)-2\sum\limits_{i=1}^N\int_0^{u_0(t_i)}e^{H(t_i)}\tilde{I}_i(t){\rm d}t.
\end{eqnarray}
By $(I_3')$ and (\ref{eq1.156}), we get
\begin{eqnarray*}
0\le \sum\limits_{i=0}^N\int_{s_i}^{t_{i+1}}e^{H(t)}(\tilde{f}_i(t,u_0(t))u_0(t)-2\tilde{F}_i(t,u_0(t))){\rm d}t,
\end{eqnarray*}
which, together with the condition $(F_2')(i)$, yields
$$
\sum\limits_{i=0}^N\int_{s_i}^{t_{i+1}}e^{H(t)}(\tilde{f}_i(t,u_0(t))u_0(t)-2\tilde{F}_i(t,u_0(t))){\rm d}t =0.
$$
Then it follows from  $(F_2')(i)$ and Remark 3.2 that $u_0\equiv 0$.
On the other hand, since $\|u_n\|_E\geq b_0>0$ for all $n$ and $u_n\to u_0$ in $E$, the continuity of the norm in $E$ ensures $\|u_0\|_E=\lim\limits_{n\to\infty}\|u_n\|_E \geq b_0>0$, which contradicts the fact that $u_0\equiv0$. Thus, no such sequences $\{u_n\}$ and $\{\varepsilon_n\}$ can exist, completing the proof.
\qed

\vskip2mm
\noindent
{\bf Remark 3.3.}  {\it By Lemma 3.6, if $u$ is a critical point of $\tilde{J}_\varepsilon(u)$ with $|\varepsilon|\leq \sigma(\frac{\sqrt{e^{H_0}}\delta}{2})$ and $\tilde{J}_\varepsilon(u)\leq \sigma(\frac{\sqrt{e^{H_0}}\delta}{2})$, then $\|u\|_E\leq \frac{\sqrt{e^{H_0}}\delta}{2}$. Furthermore, by  Lemma 2.1, we obtain that $\|u\|_{\infty}\leq\frac{\delta}{2}$. Thus, $u$ is a critical point of the original problem (\ref{eq1.1}).}

 \vskip2mm
 \noindent
{\bf Proof of Theorem 1.1.}  Without loss of generality, assume that $\varepsilon>0$. The case $\varepsilon<0$ can be addressed similarly by substituting $\tilde{g_i}(t,u)$ by $-\tilde{g_i}(t,u)$.
Next, we are prepared to verify that the functional $\tilde{J}_\varepsilon(u)$ satisfies conditions $(B_1)$-$(B_5)$ outlined in Lemma 2.3.  The conditions $(B_1)$, $(B_3)$ and $(B_5)$ follow from Lemma 3.3, Lemma 3.4 and Lemma 3.5, respectively. By $(F_1^{\prime})$ and $(I_1^{\prime})$, it is easy to see that $\tilde{J}_0(u)$ is even, so condition $(B_4)$ holds.
For $\varepsilon\in [0,1]$, by  $(F'_2)(ii)$ and (\ref{eqq3.33}), we have
\begin{eqnarray*}
|\tilde{J}_\varepsilon(u)-\tilde{J}_0(u)|\leq |\varepsilon|\sum_{i=1}^N \int_{s_i}^{t_i+1}|\tilde{G}_i(t,u(t))|{\rm d}t
\leq|\varepsilon|\sum_{i=1}^N\int_{\{t\in[s_i,t_{i+1}]:|u|\le\delta\}}|\tilde{G}_i(t,u(t))|{\rm d}t\leq|\varepsilon|C_0
:= \psi(\varepsilon),
\end{eqnarray*}
where $C_0=\sum_{i=1}^N\max_{(t,s)\in [0,T]\times[-\delta,\delta]}|\tilde{G}_i(t,s)|T$. So, condition $(B_2)$ is satisfied. Therefore, $\tilde{J}_\varepsilon(u)$ satisfies all conditions in Lemma 2.3. By using the similar arguments as Corollary 1.1 in \cite{Huang C 2022} or Theorem 1.1 in \cite{Kajikiya R 2013}, for any $\sigma>0$ and any given $k\in \mathbb{N} $, there exist $k$ distinct critical values of $\tilde{J}_\varepsilon(u)$ satisfying
\begin{eqnarray*}
-\sigma<a_{n(1)}(\varepsilon)<a_{n(2)}(\varepsilon)<\cdot\cdot\cdot<a_{n(k)}(\varepsilon)<0.
\end{eqnarray*}
Finally, owing to the arbitrariness of $\sigma$, we may choose $0<\sigma<\sigma(\frac{\sqrt{e^{H_0}}\delta}{2})$. Then by Remark 3.3, we conclude that the original problem (\ref{eq1.1}) possesses at least $k$ solutions whose $L^{\infty}$-norms are less than $\frac{\delta}{2}$. The proof is complete.
\qed

\section{Example}
  \setcounter{equation}{0}
Let $N=2$, $T=1$, and $h(t)\equiv1$ for all $t\in [0,1]$. Then $H(t)=\int_0^t 1d\tau=t$, $H_0=\min_{t\in [0,1]}{H(t)}=0$, $M_1=\max_{t\in [0,1]}{e^{H(t)}}=e$. Consider the following fourth-order ordinary differential equations:
\begin{equation}
 \label{eq11.111}
 \begin{cases}
 u^{(4)}(t)+2u^{\prime \prime \prime}(t)+\left(1+\beta e^{-t}\right)u^{\prime \prime}(t)=f(u(t)),\\
 u^{\prime \prime}(t_i^+)=u^{\prime \prime}(t_i^-),s^{\prime \prime}(t_i^+)=s^{\prime \prime}(t_i^-),u^{\prime \prime \prime}(s_i^+)=u^{\prime \prime \prime}(s_i^-),\; i=1,2,\\
 \Delta u^{\prime \prime \prime}(t_i)=I_i(u(t_i)),\; i=1,2,\\
 u(0)=u^{\prime}(0)=u(T)=u^{\prime}(T)=0,
 \end{cases}
\end{equation}
where $0=s_0<t_1=\frac{1}{5}<s_1=\frac{2}{5}<t_2=\frac{3}{5}<s_2=\frac{4}{5}<t_3=1$, $\beta<\frac{2e^{H_0}}{T^2}=2$, $f,I_1,I_2:[-\delta,\delta]\to\mathbb{R}$ with $\delta<1$ is defined by
\begin{equation}\label{4.01}
f(u)=\left\{\begin{aligned}
&f_0(u)=f_2(u)=|u|^{-\frac{1}{2}} u, \; && t\in (0,\frac{1}{5}]\cup (\frac{4}{5},1],\\
&f_1(u)=|u|^{-\frac{1}{2}} u\left(1+0.1\sin|u|\right),\; && t\in (\frac{2}{5},\frac{3}{5}],\\
&0,\; && t\in (\frac{1}{5},\frac{2}{5}]\cup(\frac{3}{5},\frac{4}{5}],\\
\end{aligned}\right.
\end{equation}
$I_1(u)=|u|^{-\frac{1}{2}} u$ and $I_2(u)=|u|^{-\frac{2}{3}} u$. Correspondingly, it is not difficult to obtain that
\begin{equation}\label{4.02}
F(u)=\int_0^{u}f(s)ds=\left\{\begin{aligned}
&F_0(u)=F_2(u)=\frac{2}{3}|u|^{\frac{3}{2}}\approx 0.667|u|^{\frac{3}{2}}, \; && t\in (0,\frac{1}{5}]\cup (\frac{4}{5},1],\\
&F_1(u)=\frac{2}{3}|u|^{\frac{3}{2}}+0.1\int_0^{u}|s|^{\frac{1}{2}}sgn(s) \sin|s|ds,\; && t\in (\frac{2}{5},\frac{3}{5}],\\
&0,\; &&  t\in (\frac{1}{5},\frac{2}{5}]\cup(\frac{3}{5},\frac{4}{5}],\\
\end{aligned}\right.
\end{equation}
From (\ref{4.02}), if $u\in [0, \delta]$, we have
\begin{eqnarray}\label{4.03}
\frac{2}{3}u^{\frac{3}{2}}+0.1\int_0^{u}s^{\frac{1}{2}} \sin sds
&=& \frac{2}{3}u^{\frac{3}{2}}+ 0.1\sin|u(\xi_1)| \int_0^{u}s^{\frac{1}{2}}ds\nonumber\\
&\geq& \frac{2}{3}u^{\frac{3}{2}}-0.1\int_0^{u}s^{\frac{1}{2}}ds\nonumber\\
&=& (\frac{2}{3}-\frac{0.2}{3})u^{\frac{3}{2}}\nonumber\\
&=& 0.6 u^{\frac{3}{2}},
\end{eqnarray}
where $u(\xi_1)\in [0,u]$. If $u\in [-\delta, 0]$, then we have
\begin{eqnarray}\label{4.04}
\frac{2}{3}(-u)^{\frac{3}{2}}+0.1\int_0^{u}(-s)^{\frac{1}{2}}(-1) (-\sin s)ds
&=& \frac{2}{3}(-u)^{\frac{3}{2}}- 0.1\int_u^{0} (-s)^{\frac{1}{2}}  \sin s ds \nonumber\\
&=& \frac{2}{3}(-u)^{\frac{3}{2}}- 0.1\sin|u(\xi_2)|\int_u^{0} (-s)^{\frac{1}{2}}ds \nonumber\\
&\geq& \frac{2}{3}(-u)^{\frac{3}{2}}- 0.1 \int_u^{0} (-s)^{\frac{1}{2}}ds \nonumber\\
&=& (\frac{2}{3}-\frac{0.2}{3})(-u)^{\frac{3}{2}}\nonumber\\
&=& 0.6 (-u)^{\frac{3}{2}},
\end{eqnarray}
where $u(\xi_2)\in [u, 0]$. Together (\ref{4.03}) with (\ref{4.04}), we obtain that
\begin{eqnarray}\label{4.05}
\frac{2}{3}|u|^{\frac{3}{2}}+0.1\int_0^{u}|s|^{\frac{1}{2}}sgn(s) \sin|s|ds\geq 0.6|u|^\frac{3}{2}.
\end{eqnarray}
Next, we verify that $f$ and $I$ satisfy those assumptions in Corollary 1.2.\\
{\bf $\bullet$} Obviously, both $f_i(u) (i=0,1,2)$ and $I_i(u) (i=1,2)$ are odd about $u\in [-\delta,\delta]$. So $(F_1)$ and $(I_1)$ hold.\\
{\bf $\bullet$} By (\ref{4.01}), (\ref{4.02}) and (\ref{4.03}), it is easy to see that $uf_i(t,u)-2F_i(t,u)\leq 1.1|u|^{\frac{3}{2}}- 1.2 |u|^{\frac{3}{2}}<0 $. Therefore, $f_i (i=0,1,2)$ satisfy $(F_2)$.\\
{\bf $\bullet$} It follows from (\ref{4.02}) and (\ref{4.05}) that condition $(F_3)$ is valid.\\
{\bf $\bullet$} If we choose $a_1=1.1$ and $\theta_1=\frac{3}{2}$, then by (\ref{4.01}), we have $|f_i(u)|\leq \max\{|u|^{\frac{1}{2}}, |u|^{\frac{1}{2}}\big(1+0.1\sin|u|\big)\}\leq a_1|u|^{\frac{1}{2}}$, thus, $(F_4)$ holds.\\
{\bf $\bullet$} If we let $a_*=0.56<\frac{a_1}{\theta_1}\approx 0.733$. By (\ref{4.02}) and (\ref{4.05}), we have $(F_5)$ holds.\\
{\bf $\bullet$} If we let $\theta_2=\frac{4}{3}$, $a_2=5$, then it is easy to see that $I_i (i=1,2)$ satisfy $(I_2)$.\\
{\bf $\bullet$} Since $I_1(u)=|u|^{-\frac{1}{2}} u\geq0$ and $I_2(u)=|u|^{-\frac{2}{3}} u\geq0$ for $u\in[0, \delta]$. Furthermore,
\begin{eqnarray*}
2\sum\limits_{i=1}^N\int_0^{u}e^{H(t_i)}I_i(t){\rm d}t &=& 2\left(\int_0^{u}e^{H(t_1)}I_1(t){\rm d}+\int_0^{u}e^{H(t_2)}I_2(t){\rm d}\right)\\
&=& 2\left(e^{H(t_1)}\frac{2}{3}|u|^{\frac{3}{2}}+e^{H(t_2)}\frac{3}{4}|u|^{\frac{4}{3}}\right)\\
&=&  \frac{4}{3}e^{\frac{1}{5}}|u|^{\frac{3}{2}}+\frac{3}{2}e^{\frac{3}{5}}|u|^{\frac{4}{3}},
\end{eqnarray*}
and
$
\sum\limits_{i=1}^Ne^{H(t_i)}I_i(u)u=e^{\frac{1}{5}}|u|^{\frac{3}{2}}+e^{\frac{3}{5}}|u|^{\frac{4}{3}},
$
so conclusion $(I_3)$ remains valid.\\
Then by Corollary 1.2, system (\ref{eq11.111}) has infinitely many nontrivial weak solutions.

\vskip2mm
\renewcommand\refname{References}
{}
\end{document}